\documentclass[letterpaper, 10 pt, conference]{ieeeconf}

\IEEEoverridecommandlockouts 
\usepackage{amsmath, amssymb, mathtools, mathrsfs, xcolor}
\usepackage[american]{circuitikz}
\usepackage{subcaption}
\usepackage{pgfplots}
\usepackage{tikz}

\ctikzset{bipoles/length=1.4cm}

\newcommand{\bbR}{\mathbb{R}}
\newcommand{\bbC}{\mathbb{C}}

\newcommand{\bbZ}{\mathbb{Z}}

\newcommand{\bbT}{\mathbb{T}}

\newcommand{\B}{\mathscr{B}}
\newcommand{\D}{\mathscr{D}}

\newcommand{\calA}{\mathcal{A}}
\newcommand{\calB}{\mathcal{B}}
\newcommand{\calC}{\mathcal{C}}
\newcommand{\calD}{\mathcal{D}}

\newcommand{\calF}{\mathcal{F}}
\newcommand{\calG}{\mathcal{G}}
\newcommand{\calH}{\mathcal{H}}

\newcommand{\calQ}{\mathcal{Q}}

\newcommand{\calX}{\mathcal{X}}

\newcommand{\calT}{\mathcal{T}}
\newcommand{\calY}{\mathcal{Y}}
\newcommand{\calS}{\mathcal{S}}
\newcommand{\calU}{\mathcal{U}}

\newcommand{\calZ}{\mathcal{Z}}

\newcommand{\frakH}{\mathfrak{H}}

\renewcommand{\epsilon}{ \varepsilon}

\newcommand{\bbm}{\begin{bmatrix*}}
\newcommand{\ebm}{\end{bmatrix*}}

\newcommand{\inv}{^{-1}}

\renewcommand{\a}{^\ast}

\newcommand{\abs}[1]{\lvert#1\rvert}
\newcommand{\norm}[1]{\lVert#1\rVert}
\newcommand{\inner}[2]{\langle#1, #2\rangle}

\renewcommand{\d}{\textrm{d}}

\newcommand{\dtk}[1]{\frac{\d^{#1}}{\d t^{#1}}}

\newcommand{\sys}{\Sigma}

\newcommand{\bbRp}{\bbR_+}
\newcommand{\bbZp}{\bbZ_{+}}
\newcommand{\intrp}{\int_{\bbR_+}}
\newcommand{\intrs}{\int_{\bbR^s}}

\newcommand{\warn}[1]{{\color{black} #1}}
\newcommand{\new}[1]{{\color{black} #1}}

\newtheorem{example}{Example}
\newtheorem{definition}{Definition}
\newtheorem{lemma}{Lemma}
\newtheorem{theorem}{Theorem}

\newtheorem{remark}{Remark}

\title{\LARGE \bfseries
	Linear time-and-space-invariant relaxation systems
}

\author{Tihol Ivanov Donchev$^1$, \quad Brayan M. Shali$^1$, \quad Rodolphe Sepulchre$^{1,2}$
	\thanks{$^\ast$The research leading to these results has received funding from the European Research Council under the Advanced ERC Grant Agreement SpikyControl n.101054323. Email:  {\itshape tiholivanov.donchev@kuleuven.be}, {\itshape brayan.shali@kuleuven.be}, {\itshape rodolphe.sepulchre@kuleuven.be}.}%
	\thanks{$^1$Department of Electrical Engineering (ESAT), KU Leuven, KasteelPark Arenberg 10, B-3001 Leuven, Belgium.
	}
	\thanks{$^2$Department of Engineering, University of Cambridge, TrumpingtonStreet, Cambridge CB2 1PZ, United Kingdom.}%
}

\begin{document}
	\maketitle
	\thispagestyle{empty}
	\pagestyle{empty}
	
	\begin{abstract}
		This paper generalizes the physical property of relaxation from linear time-invariant (LTI) to linear time-and-space-invariant (LTSI) systems. It is shown that the defining features of relaxation---complete monotonicity, passivity, and memory-based storage---carry over seamlessly to the spatio-temporal domain. An LTSI system is shown to be of relaxation type if and only if its associated spatio-temporal Hankel operator is cyclically monotone. This implies the existence of an intrinsic quadratic storage functional defined uniquely by past inputs, independently of any state-space realization. As in the LTI case, LTSI relaxation systems are shown to be those systems for which the state-space concept of storage coincides with the input-output concept of \warn{fading} memory functional.
		
	\end{abstract}

	\section{Introduction}
	
	Relaxation systems have recently received renewed attention in systems and control theory due to their relevance in structured optimal control, circuit synthesis, and modeling of physical and biological processes \cite{pates2019, pates2022, chaffey2023, rodolphe2024, vanderschaft2024}. Yet, existing literature mostly addresses relaxation from a purely \textit{temporal} perspective. Many real-world applications---particularly in distributed-parameter systems governed by partial differential equations---require a framework that also accounts for \textit{spatial} variations. In this paper, we propose such an extension, generalizing the classical theory of linear time-invariant (LTI) relaxation systems to systems that depend on both time and space. Specifically, we address systems that are also spatially invariant, resulting in a linear time-and-space-invariant (LTSI) system.
	
	In the LTI setting, the relaxation property admits two classical characterizations \cite{willems1972b}. From a state-space viewpoint, relaxation systems are (internally) reciprocal and passive, containing only one type of energy storage element. From an input-output perspective, relaxation systems correspond to convolution operators whose impulse responses are completely monotone, i.e., their memory decays monotonically and without a hint of oscillation. The recent work \cite{chaffey2023} shows that relaxation can be characterized as cyclic monotonicity of the Hankel operator. The cyclic monotonicity of the Hankel operator implies the existence of an intrinsic storage functional uniquely determined  by past inputs alone. The latter is consistent with Willems' notion of storage, hence this result identifies relaxation systems as those systems for which the state-space notion of storage coincides with the input-output notion of memory.
	
	More recently, these ideas have been extended to the nonlinear setting. Building on the state-space notions of reciprocity and passivity, \cite{vanderschaft2024} defines nonlinear relaxation systems as gradient systems with a Hessian Riemannian metric satisfying both reciprocity and passivity conditions, thus ensuring energy dissipation without oscillatory behavior. An alternative definition is provided in \cite{rodolphe2024}, where nonlinear relaxation systems are characterized from an input-output perspective as fading memory operators with completely monotone impulse responses. This characterization relies on the universal approximation property of nonlinear fading memory operators \cite{boyd1985}, which shows that a nonlinear fading memory operator can be uniformly approximated by an LTI system composed with a nonlinear static readout. The input-output characterization introduced in \cite{rodolphe2024} is consistent with the state-space definition of \cite{vanderschaft2024} in the sense that the former can be seen as a special case of the latter.
	
	In this paper, we show that the defining features of relaxation---complete monotonicity, passivity, and memory-based storage---carry over seamlessly to the spatio-temporal domain. Specifically, we prove that an LTSI system is of relaxation type precisely when its spatio-temporal Hankel operator is cyclically monotone. As such, the latter is the gradient of a closed convex functional defined uniquely by past inputs alone. This establishes a direct link between the internal (state-space) and external (input-output) perspectives on energy storage. The results are illustrated throughout this paper with the diffusion equation as a canonical example of spatio-temporal relaxation.

	The remainder of this paper is organized as follows. In Section~\ref{sec:diffusion_relaxation}, we \new{introduce} diffusion as a natural example of spatio-temporal relaxation. Section~\ref{sec:notation_preliminaries} contains notation and preliminaries. In Section~\ref{sec:ltsi_relaxation}, we define and characterize (internal) relaxation for LTSI systems.  In Section~\ref{sec:passivity_relaxation}, we show that internal relaxation implies impedance passivity with storage determined from the external behaviour. In Section VI, we formalize the connection between internal and external by establishing properties of the Hankel operator of an LTSI relaxation system. Finally, Section~\ref{sec:conclusion} summarizes our findings and discusses directions for further work.

	\section{Diffusion and relaxation}\label{sec:diffusion_relaxation}
	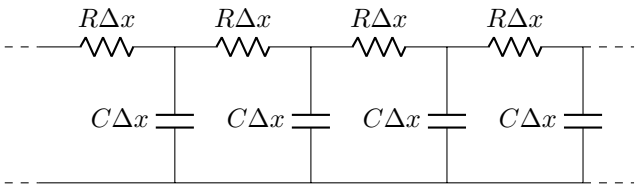
\begin{figure}[b]
		\centering
		\ctikzset{resistors/scale=0.6, inductors/scale=0.6, capacitors/scale=0.6}
		
		\begin{circuitikz}[scale=0.9]
			\draw[dashed]
			(-.5,2) to (0,2)
			(-.5,0) to (0,0)
			(8,2) to (8.75,2)
			(8,0) to (8.75,0);
			\draw
			(0,2) to[R=$R\Delta x$] (2,2)
			(2,2) to[R=$R\Delta x$] (4,2)
			(4,2) to[R=$R\Delta x$] (6,2)
			(6,2) to[R=$R\Delta x$] (8,2)
			
			(0,0) to (8,0)
			
			(2,0) to[C=$C\Delta x$, left] (2,2)
			(4,0) to[C=$C\Delta x$, left] (4,2)
			(6,0) to[C=$C\Delta x$, left] (6,2)
			(8,0) to[C=$C\Delta x$, left] (8,2);
			
		\end{circuitikz}
		\caption{\small Distributed RC line circuit.}
		\label{fig:RCdistrib}
		\vspace{-3mm}
	\end{figure}
	
\begin{figure}[h!]
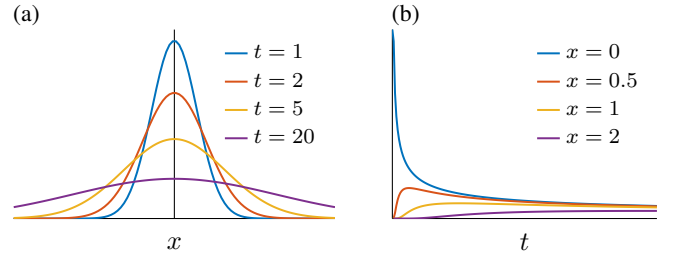

		\captionsetup[subfigure]{skip=1pt, margin=0 mm, singlelinecheck=false}
		\centering
		\begin{subfigure}[t]{0.58\linewidth}
			\subcaption{}
			\input{heat_kernel_fixed_t}
			\label{fig:diffusion_plots_a}
		\end{subfigure}%
		\begin{subfigure}[t]{0.5\linewidth}
			\subcaption{}
			\input{heat_kernel_fixed_x}
			\label{fig:diffusion_plots_b}
		\end{subfigure}
		\caption{\small Plots of the heat kernel \eqref{eq:heat_kernel} at various times (a) and various locations (b).}
		\label{fig:diffusion_plots}
		\vspace{-3mm}
	\end{figure}

	\new{To motivate the content of this paper, we argue that diffusion is the physical phenomenon that should underline the system theoretic property of relaxation for LTSI systems.} Historically, relaxation was first \new{identified \cite{meixner1964} as the system theoretic property of passive electrical circuits with only one type of energy storage element, i.e., RC or RL circuits.} \new{This characterization naturally extends} to spatially distributed systems by considering an infinite transmission line composed of resistors and capacitors periodically arranged along a spatial axis.
	Figure~\ref{fig:RCdistrib} illustrates the schematic of such an infinite distributed RC transmission line. Taking the spatial discretization interval $\Delta x \rightarrow 0$, this discrete electrical network is precisely governed by the one-dimensional diffusion equation \cite{kaufman1962}, i.e. 
	\begin{equation}\label{eq:diffRC}
		C\frac{\partial}{\partial t}V(t,x) = \frac{1}{R}\frac{\partial^2}{\partial x^2}V(t,x) + I(t,x),
	\end{equation}
	where $V(t,x)$ and $I(t,x)$ \new{are the voltage and current} at time $t$ and position $x$. \new{Denoting} the variables $y(t,x)=V(t,x)$, $u(t,x)=I(t,x)/C$, and the coefficient $\alpha=1/RC$, the latter simplifies to the standard (heat) diffusion equation
	\begin{equation}\label{eq:diff}
		\frac{\partial}{\partial t}y(t,x) = \alpha\frac{\partial^2}{\partial x^2}y(t,x) + u(t,x).
	\end{equation}
	The explicit solution of \eqref{eq:diff}, see, e.g., \cite{evans2010}, is given by the convolution
	\begin{equation}\label{eq:diff_solution}
		y(t,x) = \int_{-\infty}^t \int_{-\infty}^\infty g(t-\tau,x-\xi)u(\tau,\xi)\, \d \xi \d \tau,
	\end{equation}
	where $g$ is the so-called \emph{heat kernel} given by
	\begin{equation}\label{eq:heat_kernel}
		g(t,x) = \frac{1}{\sqrt{4\pi \alpha t}}e^{-\frac{x^2}{4\alpha t}}.
	\end{equation}
	The heat kernel clearly exhibits the relaxation property \new{when regarded as a spatially extended function of time.} In particular, Figure~\ref{fig:diffusion_plots_a} illustrates the spatial profile of \eqref{eq:heat_kernel} at various times $t$. It shows the relaxation property, as the spatial profile flattens monotonically towards equilibrium as $t$ increases. However, Figure~\ref{fig:diffusion_plots_b} shows the kernel’s temporal behavior at a fixed location, highlighting that the kernel is not a completely monotone function in time, as one might expect from a classical LTI relaxation system. 
	
	This apparent discrepancy is resolved by examining the spatial Fourier transform of the diffusion equation \eqref{eq:diff}, which yields the infinite family of decoupled LTI systems
	\begin{equation}
		\frac{\partial}{\partial t}\hat y(t,\omega) = -\alpha\omega^2 \hat y(t,\omega) + \hat u(t,\omega)
	\end{equation}
	parameterized by the frequency variable $\omega$. The impulse response of the latter is given by 
	\begin{equation}\label{eq:fourier_kernel} \hat g(t,\omega) = e^{-\alpha \omega^2 t}.
	\end{equation}	
	In this frequency domain representation, each mode $\omega$ individually defines an LTI system whose impulse response is completely monotone in time. Thus, the diffusion equation defines an infinite-dimensional system that can be decomposed into an infinite family of finite-dimensional LTI relaxation systems.
	
	The diffusion process provides a clear physical connection from classical relaxation theory to its infinite-dimensional, spatially extended counterpart. This is further reinforced by the recent work \cite{drummond2024}, where it is explicitly demonstrated that the spatial discretization of the diffusion equation leads to a finite-dimensional LTI relaxation system. In this paper, we extend the theory of LTI relaxation systems to a class of LTSI system, with the diffusion process as a canonical example of an LTSI relaxation system.
	
	\section{Notation and preliminaries}\label{sec:notation_preliminaries}
	
	We use mostly standard notation. We denote the sets of nonnegative real numbers and integers by $\bbRp$ and $\bbZp$, respectively. We denote the inner-product on a (complex) Hilbert space $\calX$ by $\inner{\cdot}{\cdot}_\calX$ and the induced norm by $\norm{\cdot}_\calX$. We omit the subscript if $\calX = \bbC^n$ or there is no risk of ambiguity. We denote the Banach space of bounded linear operators from $\calX$ to another Hilbert space $\calY$ by $\B(\calX,\calY)$, and the operator norm by $\norm{\cdot}$. We use the shorthand notation $\B(\calX,\calX) = \B(\calX)$. We also consider unbounded linear operators $\calA:\D(\calA) \to \calY$, where $\D(\calA) \subset\calX$ is the domain of $\calA$. We say that $\calA$ is densely defined if $\D(\calA)\subset\calX$ is dense. If $\calA$ is densely defined, then it has an adjoint $\calA\a: \D(\calA\a) \to \calX$, where $\calD(\calA\a)\subset\calY$ is closed. We say that $\calA$ is self-adjoint if $\D(\calA) = \D(\calA\a)$ and $\inner{\calA x}{y} = \inner{x}{\calA y}$ for all $x,y\in\D(\calA)$. We say that a self-adjoint $\calA$ is nonnegative if $\inner{\calA x}{x} \geq 0$ for all $x\in\D(\calA)$, and we write $\calA \geq 0$ to indicate that $\calA$ is self-adjoint and nonnegative. 
	
	\subsection{Square integrable functions and the Fourier transform}
	We denote the Hilbert space of square integrable functions from a domain $X$ to a Hilbert space $\calX$ by $L_2(X, \calX)$, where
	\begin{equation}
		\inner{f}{g}_{L_2(X,\calX)} = \int_{X} \inner{f(x)}{g(x)}_\calX\, \d x.
	\end{equation}
	The Fourier transform of a function $f:\bbR^s\to\bbC^n$ is the function $\hat f: \bbR^s\to\bbC^n$ given by
	\begin{equation}
		\hat f(\omega) = \new{(2\pi)^{-s/2}}\intrs e^{-j\inner{\omega}{x}}f(x)\, \d x.
	\end{equation}
	It is well-known \cite{zuazo2001} that the Fourier transform converges for absolutely integrable functions. Due to Plancherel's theorem, the Fourier transform can be extended to a \emph{unitary} operator $\calF\in \B(L_2(\bbR^s, \bbC^n))$, i.e., $\calF\inv = \calF\a$. Since $\calF$ agrees with the Fourier transform on the dense subspace of absolutely integrable functions, we refer to $\calF f$ as the Fourier transform of $f$, and we write $\hat f = \calF f$. 
	
	A key property of the Fourier transform that we exploit in this paper is that translation-invariant operators are transformed into multiplication operators. The Fourier transform of a densely defined operator $\calA$ on $L_2(\bbR^s,\bbC^n)$ is the densely defined operator $\hat \calA = \calF \calA \calF\a$ on $L_2(\bbR^s,\bbC^n)$, such that $\hat \calA \hat f\in L_2(\bbR^s,\bbC^n)$ is the Fourier transform of $\calA f \in L_2(\bbR^s,\bbC^n)$. For each $\bar x\in\bbR^s$, the \emph{translation operator} $\calT_{\bar x}\in\B(L_2(\bbR^s,\bbC^n))$ is given by $(\calT_{\bar x} f)(x) = f(x+\bar x)$. A densely defined operator $\calA$ is \emph{translation-invariant} if $\calT_{\bar x}:\D(\calA)\to\D(\calA)$ and $\calT_{\bar x}\calA = \calA \calT_{\bar x}$ for all $\bar x\in\bbR^s$. The Fourier transform of a translation-invariant operator is a \emph{multiplication operator} \cite[Section~3.6.7]{zuazo2001}, i.e., there exists a function $\hat A:\bbR^s\to \bbC^{n\times n}$ such that $(\hat\calA \hat f)(\omega) = \hat A(\omega) \hat f(\omega)$. The function $\hat A$ is called the \emph{symbol} of the operator $\hat\calA$. It is well-known that the operators $\hat \calA$ and, thus, $\calA$ are bounded if and only if the symbol $\hat A$ is bounded.
	
	Throughout this paper, we consider functions that depend on a temporal variable $t$ and a spatial (frequency) variable $x$ ($\omega$). We denote the dependence on the spatial (frequency) variable with a subscript, e.g., we write $f_x$ instead of $f(x)$. In particular, we often consider functions of the form $f:\bbT\to  L_2(\bbR^s,\bbC^m)$, where $\bbT$ is a temporal domain and $\bbR^s$ is a spatial (frequency) domain. We interpret such functions as mapping the temporal variable $t\in\bbT$ to a spatial (frequency) profile $f(t)\in L_2(\bbR^s,\bbC^m)$. We denote the dependence on the spatial (frequency) variable with a subscript, i.e., we write $f_x(t) = f(t)(x)$, so that the function $f_x:\bbT\to \bbC^n$ is the temporal trajectory at a given point in space $x\in\bbR^s$. Finally, we write $\hat f$ to denote the spatial Fourier transform of $f$, i.e., the function $\hat f: \bbT\to  L_2(\bbR^s,\bbC^m)$ such that $\hat f(t)\in L_2(\bbR^s,\bbC^m)$ is the Fourier transform of $f(t)\in L_2(\bbR^s,\bbC^m)$ for all $t\in\bbT$.
	
	\subsection{LTI relaxation systems}
	
	Consider the linear time-invariant (LTI) system
	\begin{equation}\label{eq:lti_ss}
		\sys:\left\lbrace \begin{aligned}
			\frac{\d }{\d t}z(t) &= Az(t) + Bu(t),\\
			y(t) &= Cz(t),
		\end{aligned}\right.
	\end{equation}
	with input space $\bbR^m$, state space $\bbR^n$, and output space $\bbR^p$. The impulse response of the system $\sys$ is the matrix-valued map $g:\bbRp\to\bbR^{p\times m}$ given by $g(t) = Ce^{At} B$. Assuming that the system $\sys$ is initially at rest, i.e., $z(-\infty) = 0$, its output is given by the convolution
	\begin{equation}
		y(t) = \int_{-\infty}^{t} g(t-\tau)u(\tau)\, \d \tau.
	\end{equation}
	The system $\sys$ is of \emph{relaxation type} \cite{willems1972b} if $m = p$ and its impulse response $g$ is devoid of any oscillatory behaviour, expressed as $g$ being \emph{completely monotone} \new{\cite[Chapter~IV]{widder1941}}, i.e.,
	\begin{equation}
		(-1)^k\dtk{k} g(t) \geq 0
	\end{equation}
	for all $t\in\bbRp$ and $k\in\bbZp$. By Bernstein's theorem \new{\cite{bernstein29}}, $g$ is completely monotone if and only if it is a conic combination of decaying exponentials, i.e., $g(t) = \sum_{i=1}^{n} G_ie^{-p_it}$, where $G_i\geq 0$ and $p_i \geq 0$ for all $i\in\{1,\dots,n\}$. Relaxation systems are passive and reciprocal with only one type of energy storage. They admit a unique storage function that can be determined from the input-output behaviour of the system. Furthermore, every relaxation system has a state space realization that is internally of relaxation type \cite{willems1976}, i.e., $B= C\a$ and $A\leq 0$. 
	
	For an exponentially stable system $\sys$, relaxation can be characterized via the Hankel operator, i.e., the bounded linear operator $\calH: L_2(\bbRp, \bbR^m)\to L_2(\bbRp, \bbR^p)$ given by
	\begin{equation}\label{eq:H_lti}
		(\calH v)(t) = \intrp g(t+\tau)v(\tau)\, \d\tau.
	\end{equation}
	The Hankel operator $\calH$ maps past inputs to future outputs in the sense that the input $u(t) = v(-t)$ for all $t\leq 0$ and $u(t) = 0$ for all $t >0$, results in the output $y(t) = (\calH v)(t)$ for all $t\geq 0$. In \cite{chaffey2023}, it is shown that the Hankel operator of a relaxation system is self-adjoint and nonnegative, hence cyclically monotone. This implies that $\calH$ is the functional derivative of a quadratic (memory) functional, which serves as an intrinsic energy storage for the system. 
	
	Finally, we note that the results in this section are also valid for LTI systems with complex input, state and output spaces, which are considered throughout this paper.
	
	\section{LTSI systems and relaxation}\label{sec:ltsi_relaxation}
	
	In this section, we define and characterize the notion of relaxation for a class of linear time-and-space-invariant (LTSI) systems. The definition of relaxation for LTSI systems is a natural generalization of the definition for LTI systems, namely, it is expressed as complete monotonicity of the (operator-valued) impulse response. The spatial invariance of the LTSI system allows us to view it as a family of LTI systems via the Fourier transform. This allows us to characterize relaxation of the LTSI system as relaxation of the corresponding family of LTI systems.
	
	Consider the LTSI system $\sys$ of the form
	\begin{equation}\label{eq:ltsi_ss}
		\sys:\left\lbrace\begin{aligned}
			\frac{\d}{\d t}z(t) &= \calA z(t) + \calB u(t),\\
			y(t) &= \calC z(t),
		\end{aligned}\right.
	\end{equation}
	with (infinite-dimensional) input space $\calU = L_2(\bbR^s, \bbC^m)$, state space $\calZ = L_2(\bbR^s, \bbC^n)$, output space $\calY = L_2(\bbR^s, \bbC^p)$, and translation-invariant linear operators $\calA:\calD(\calA) \to \calZ$, $\calB \in \B(\calU,\calZ)$ and $\calC\in\B(\calZ,\calY)$, where we assume that the domain $\D(\calA)\subset\calZ$ is dense. The system is time-invariant because the operators $\calA$, $\calB$ and $\calC$ are constant, and space-invariant because they are translation-invariant. The operator $\calA$ is typically unbounded and the notion of solution requires some care, see, e.g. \cite{curtain2020}. Without going into details, we assume that $\calA$ is the infinitesimal generator of a $C_0$-semigroup $t\mapsto \exp(\calA t)\in\B(\calZ)$. The notation is motivated by the fact that the latter is a generalization of the exponential function. 
	
	The impulse response of the LTSI system $\sys$ is the operator-valued map $\calG:\bbRp\to\B(\calU,\calY)$ given by
	\begin{equation}
		\calG(t) = \calC\new{\exp(\calA t)}\calB.
	\end{equation}
	The output of $\sys$ is given by the convolution
	\begin{equation}
		y(t) = \int_{-\infty}^{t} \calG(t-\tau)u(\tau)\, \d\tau,
	\end{equation}
	where we have assumed that $z(-\infty) = 0$, i.e., the system is initially at rest. Since $\calA$, $\calB$ and $\calC$ are translation-invariant, their Fourier transforms $\hat \calA$, $\hat\calB$ and $\hat\calC$ are multiplication operators with symbols $\hat A$, $\hat B$ and $\hat C$. Consequently, the Fourier transform of \eqref{eq:ltsi_ss} yields the family of LTI systems
	\begin{equation}
		\hat\sys_\omega: \left\lbrace
		\begin{aligned}
			\frac{\d}{\d t}\hat z_\omega(t) &= \hat A_\omega\hat z_\omega (t) + \hat B_\omega \hat u_\omega(t), \\ \hat y_\omega(t)&= \hat C_\omega\hat z_\omega(t),
		\end{aligned}\right.
	\end{equation}
	parametrized by the frequency variable $\omega\in\bbR^s$. The Fourier transform $\hat \calG(t)$ of $\calG(t)$ is a multiplication operator with symbol $\hat g(t)$ given by $	\hat g_\omega(t) = \hat C_\omega e^{\hat A_\omega t}\hat B_\omega$, where $\hat g_\omega$ is the impulse response of the LTI system $\hat\sys_\omega$. We assume that the symbols $\hat A$, $\hat B$ and $\hat C$ are continuous, so that $\hat g(t)$ is continuous for all $t\in\bbR_+$. 
	
	The fact that an LTSI system can be viewed as a family of LTI systems greatly simplifies its analysis. Indeed, many system theoretic properties of $\sys$ can be verified \emph{pointwise} by verifying the analogous property for all $\hat\sys_\omega$, $\omega\in\bbR^s$, see, e.g., \cite{bamieh2002, arbelaiz2024}. This is the overarching theme in this paper. 
	With this in mind, as a natural generalization of the definition of relaxation for LTI systems, we  define relaxation for LTSI systems as follows.
	\begin{definition}
		The LTSI system $\sys$ is of relaxation type if $\calU = \calY$ and the impulse response $\calG$ is completely monotone.
	\end{definition}
	
	The impulse response $\calG$ is completely monotone if its $k$-th strong derivative satisfies
	\begin{equation}\label{eq:impulse_cm}
		(-1)^k\frac{\d^k}{\d t^k}\calG (t) \geq 0
	\end{equation}
	for all $t >0$ and $k\in\bbR_+$. Since the Fourier transform is unitary, completely monotonicity of $\calG$ is equivalent to complete monotonicity of $\hat \calG$. The latter is a multiplication operator, hence its properties can be inferred from its symbol, i.e., the impulse responses $\hat g_\omega$, $\omega\in\bbR^s$. Consequently, relaxation can be verified pointwise, as shown in the following theorem.
	
	\begin{theorem}\label{thm:relaxation_fourier}
		The LTSI system $\sys$ is of relaxation type if and only if the LTI system $\hat \sys_\omega$ is of relaxation type for all $\omega\in\bbR^s$.
	\end{theorem}
	\proof
		As shown in \cite[Theorem~2.1.13]{curtain2020}, the semigroup $\exp(\calA t)$ is strongly differentiable and the $k$-th strong derivative is given by $\exp(\calA t)\calA^k$ for all $t > 0$ on the dense domain $\D(\calA^k)\subset\calZ$ . Therefore, the impulse response $\calG(t)$ is strongly differentiable with $k$-th strong derivative given by
		\begin{equation}
				\frac{\d^k}{\d t^k} \calG(t) = \calC\exp(\calA t)\calA^k\calB
			\end{equation}
		The Fourier transform of the latter is a multiplication  operator $\hat \calG^k(t)$ with \emph{continuous} symbol $\hat g^k(t)$ given by
		\begin{equation}
				\hat g^k_\omega(t) = \frac{\d^k}{\d t^k}\hat g_\omega(t) = \hat C_\omega e^{\hat A_\omega t}\hat A_\omega ^k\hat B_\omega
			\end{equation}
		Therefore, $(-1)^k\hat \calG^k(t) \geq 0$ if and only if $(-1)^k \hat g^k_\omega(t)\geq 0$ for all $\omega\in\bbR^s$. Moreover, since  the Fourier transform is unitary, $(-1)^k\hat \calG^k(t) \geq 0$ if and only if \eqref{eq:impulse_cm} holds. Combining the latter two statements, it follows that $\calG(t)$ is completely monotone if and only if $\hat g_\omega$ is completely monotone for all $\omega\in\bbR^s$, which concludes the proof.
	\endproof
	
	Recall that LTI relaxation systems have state-space realizations that are internally of relaxation type \cite{willems1976}. This is not necessarily true for LTSI relaxation systems, mainly due to the technical assumptions on the operators that define an LTSI system, see Remark~\ref{rem:relaxation_internal_relaxation}. Whenever necessary, we avoid this issue by simply assuming that $\sys$ is internally of relaxation type, defined below.
	\begin{definition}
		The LTSI system $\sys$ is \emph{internally of relaxation type} if  $\calA \leq 0$ and $\calB = \calC\a$.
	\end{definition}
	
	Analogously to Theorem~\ref{thm:relaxation_fourier}, we obtain the following characterization of internal relaxation.
	\begin{theorem}\label{thm:internal_relaxation_fourier}
		The LTSI system $\sys$ is internally of relaxation type if and only if the LTI system $\hat \sys_\omega$ is internally of relaxation type for all $\omega\in\bbR^s$.
	\end{theorem}
	\proof
		Since the Fourier transform is unitary, $\calA \leq 0$  and $\calB=\calC\a$ if and only if $\hat\calA \leq 0$ and $\hat\calB =\hat\calC\a$. Since the symbols $\hat A$, $\hat B$, and $\hat C$ of the multiplication operators $\hat\calA$, $\hat \calB$ and $\hat\calC$ are continuous, $\hat\calA \leq 0$ and $\hat\calB =\hat\calC\a$ if and only if $\hat A_\omega \leq 0$ and $\hat B_\omega = \hat C_\omega\a$ for all $\omega\in\bbR^s$.
	\endproof
	
	We now go back to the diffusion equation as an example, after which we conclude this section with a couple of remarks.
	\begin{example}\label{ex:diff_relaxation}
		As expected, diffusion is an example of LTSI relaxation. Indeed, note that the diffusion equation \eqref{eq:diff} is an LTSI system of the form \eqref{eq:ltsi_ss} with $\calA = \alpha\frac{\partial^2}{\partial x^2}$ and $\calC=\calB$ given by the identity operator. We have already seen that the impulse response \eqref{eq:fourier_kernel} is completely monotone for all $\omega\in\bbR^s$, hence the diffusion equation is an LTSI system of relaxation type due to Theorem~\ref{thm:relaxation_fourier}. Furthermore, it is well-known that the second-order differential operator $\frac{\partial^2}{\partial x^2}$ is self-adjoint and nonpositive, hence the diffusion equation is an example of an LTSI system that is \emph{internally} of relaxation type.
	\end{example}
	
	\begin{remark}\label{rem:relaxation_internal_relaxation}
		It is easily seen that internal relaxation implies relaxation. For the converse, suppose that the LTSI system $\sys$ is of relaxation type, hence, due to Theorem~\ref{thm:relaxation_fourier}, the LTI system $\hat \sys_\omega$ is of relaxation type for all $\omega\in\bbR^s$. This means that $\hat \sys_\omega$ has a state-space realization that is internally of relaxation type, i.e., $\hat g_\omega(t) = {\bar{B}}_\omega \a e^{{\bar A}_\omega t} {\bar{B}}_\omega$, where ${\bar A}_\omega \leq 0$. Using this, we can define $\bar\calA$ and $\bar \calB$ as the translation-invariant operators whose Fourier transforms are the multiplicative operators with symbols ${\bar A}$ and ${\bar B}$. Then, we can write $\calG(t) = \bar \calB\a \exp(\bar\calA t) \bar\calB$ so that $\calG$ can be seen as the impulse response of an LTSI system that is internally of relaxation type. The issue here is that $\bar\calA$ is not guaranteed to be the infinitesimal generator of a $C_0$-semigroup, and $\bar \calB$ is not guaranteed to be bounded. Nevertheless, examples where these conditions are not satisfied seem to be artificial mathematical constructs rather than real physical examples. Similar issues arise when dealing with (impedance) passivity, see Remark~\ref{rem:relaxation_passivity}.
	\end{remark}
	
	\begin{remark}
		Systems of the form \eqref{eq:ltsi_ss} with $\calB = \calC\a$ are known as \emph{collocated systems} \cite[p.266]{curtain2020}. The name reflects the practice of applying control and observation actions at the same point in distributed parameter systems. We refer to a collocated system of the form \eqref{eq:ltsi_ss} with $\calA = \calA\a$ as \emph{internally symmetric}. If, in addition, $\calA \leq 0$, then $\calA$ is the infinitesimal generator of a \emph{contraction semigroup}, see \cite[Section~2.3]{curtain2020} for details. The latter is a mild form of stability which is implied by, e.g., exponential stability. Therefore, internal relaxation in LTSI systems is a combination of internal symmetry and stability, just like internal relaxation in LTI systems \cite{willems1972b, willems1976}. 
	\end{remark}
	
	\section{Impedance passivity and relaxation}\label{sec:passivity_relaxation}
	
	In this section, we show that LTSI systems that are internally of relaxation type are impedance passive with energy storage that can be determined from their input-output behaviour, thus extending the analogous result on LTI relaxation systems \cite{willems1972b}. In doing so, we also characterize impedance passivity for LTSI systems, namely, we show that it can be verified pointwise, just like (internal) relaxation. 
	
	To begin with, impedance passivity is the analogue of passivity for infinite-dimensional systems, see \cite[Section~7.5]{curtain2020} for details. It is defined as follows.
	\begin{definition}
		The LTSI system $\sys$ is \emph{impedance passive} if $\calU = \calY$ and there exists a storage functional $\calS:\calZ\to\bbR_+$, given by $\calS(z) = \inner{\calQ z}{z}$ with $\calQ\in\B(\calZ)$, $\calQ\geq 0$, such that 
		\begin{equation*}
			\calS(z(t)) \leq \calS(z(0)) + \int_{0}^{t}  \inner{u(\tau)}{y(\tau)} + \inner{y(\tau)}{u(\tau)}\, \d t
		\end{equation*}
		for all $z(0)\in\calZ$, $t\geq 0$, and $u\in L_2([0,\tau], \calU)$.
	\end{definition}
	
	Impedance passivity of infinite-dimensional systems has a similar characterization as passivity of finite-dimensional systems, namely, it can be expressed as the solvability of a linear operator inequality.
	\begin{lemma}\label{lem:passivity_ltsi}
		The LTSI system $\sys$ is impedance passive if and only if there exists $\calQ\in\B(\calZ)$, $\calQ\geq 0$, such that 
		\begin{equation}\label{eq:passivity_ltsi}
			\calC = \calB\a \calQ, \qquad \inner{\calA z}{\calQ z} + \inner{\calQ z}{\calA z} \leq 0
		\end{equation}
		for all $z\in\calD(\calA)$.
	\end{lemma}
	\proof See \new{\cite[Lemma~7.5.4]{curtain2020}}. \endproof
	
	Just like (internal) relaxation, impedance passivity of LTSI systems can be verified pointwise. This is shown in the following theorem.
	\begin{theorem}
		The LTSI system $\sys$ is impedance passive if and only if there exist $\hat Q_\omega\in\bbC^{n\times n}$, $\hat Q_\omega \geq 0$, such that
		\begin{equation}\label{eq:passivity_omega}
			\hat C_\omega = \hat B_\omega\a \hat Q_\omega,, \qquad \hat A_\omega\a \hat Q_\omega + \hat Q_\omega \hat A_\omega \leq 0, 
		\end{equation} 
		for all $\omega\in\bbR^s$, and $\sup_{\omega\in\bbR^s} \norm{\hat Q_\omega} < \infty$.
	\end{theorem}
	\proof
		This proof is similar to the proof of \cite[Theorem~2]{bamieh2002}. The idea is that the condition $\sup_{\omega\in\bbR^s} \norm{\hat Q_\omega} < \infty$ is equivalent to the boundedness of the multiplication operator with symbol $\hat Q$, which we denote by $\hat \calQ\in\B(\calZ)$. We start by proving sufficiency. Since $\hat Q_\omega \geq 0$ and \eqref{eq:passivity_omega} hold for all $\omega\in\bbR^s$, it follows that $\hat \calQ \geq 0$ and
		\begin{equation}\label{eq:hatQ}
				\hat \calC = \hat \calB\a \hat \calQ,\qquad \inner{\hat\calA \hat z}{\hat\calQ \hat z} + \inner{\hat\calQ \hat z}{\hat\calA \hat z} \leq 0
			\end{equation}
		for all $\hat z\in\calD(\hat \calA)$. Recall that the Fourier transform is unitary. This implies that $\calQ = \calF\a \hat\calQ \calF \in \B(\calZ)$, satisfies $\calQ \geq 0$ and the conditions of Lemma~\ref{lem:passivity_ltsi}, hence $\sys$ is impedance passive. 
		
		For the proof of necessity, suppose that $\sys$ is impedance passive. \warn{Since the dynamics of $\sys$ are translation-invariant, the storage is translation-invariant, i.e., there exists translation-invariant} $\calQ\in\B(\calZ), \calQ \geq 0,$ that satisfies the conditions of Lemma~\ref{lem:passivity_ltsi}.
		This implies that $\hat\calQ = \calF \calQ \calF\a \in\B(\calZ)$ is a multiplication operator that satisfies $\hat \calQ\geq 0$ and \eqref{eq:hatQ}. Finally, the symbol $\hat Q$ of $\hat \calQ$ is such that $\sup_{\omega\in\bbR^s} \norm{\hat Q_\omega} < \infty$  because $\hat\calQ\in\B(\calZ)$, $\hat Q_\omega \geq 0$ for all $\omega\in\bbR^s$ because $\hat\calQ \geq 0$, and \eqref{eq:passivity_omega} holds for all $\omega\in\bbR^s$ because \eqref{eq:hatQ} holds.
	\endproof
	
	In \cite{willems1972b}, it is shown that LTI relaxation systems are passive with a unique\footnote{assuming minimality of the state-space realization} compatible storage function that can be determined from their input-output behaviour. Deriving the analogue of this result for LTSI relaxation systems presents some minor technical challenges, see Remark~\ref{rem:relaxation_passivity}. We circumvent these by instead considering LTSI systems that are \emph{internally} of relaxation type. In particular, we obtain the following theorem.
	\begin{theorem}\label{thm:storage_io}
		If the LTSI system $\sys$ is internally of relaxation type, then it is impedance passive with storage $\calS:\calZ\to\bbRp$ given by $\calS(z) = \norm{z}^2$. Furthermore,
		\begin{equation}\label{eq:ltsi_relaxation_passive}
			\calS(z(0))  = \int_{0}^{\infty} \inner{u(-t)}{y(t)}\, \d t,
		\end{equation}
		where $z:\bbR\to\calZ$ and $y:\bbR\to\calY$ are the state and output trajectories corresponding to the input trajectory $u:\bbR\to\calU$ such that $u(t) = 0$ for $t \geq 0$.
	\end{theorem}
	\proof
		Since $\sys$ is internally of relaxation type, we have that $\calA \leq 0$ and $\calB = \calC\a$. Let $\calQ\in\B(\calZ)$ be the identity operator and note that $\calQ \geq 0$ and $\calS(z) = \inner{\calQ z}{z}$. Furthermore, \eqref{eq:passivity_ltsi} holds for all $z\in\calZ$, hence $\sys$ is impedance passive due to Lemma~\ref{lem:passivity_ltsi}. Now, consider an input trajectory $u:\bbR\to\calU$ such that $u(t) = 0$ for $t \geq 0$. Let $z:\bbR\to\calZ$ and $y:\bbR\to\calY$ be the resulting state and output trajectories. Consider the Fourier transforms $\hat u$, $\hat z$ and $\hat y$ of $u$, $z$ and $y$, respectively. Note that $\hat z_\omega$ and $\hat y_\omega$ are the state and output trajectories of $\hat \sys_\omega$ corresponding to the input trajectory $\hat u_\omega$, which is such that $\hat u_\omega(t) = 0$ for $t\geq 0$. Due to Theorem~\ref{thm:internal_relaxation_fourier}, $\hat \sys_\omega$ is internally of relaxation type for all $\omega\in\bbR^s$ and, thus, $\hat A_\omega \leq 0$ and $\hat B_\omega = \hat C_\omega\a$. As shown in \cite[Theorem~9]{willems1972b}, see also \cite[Remark~8]{willems1972b}, this implies that 
		\begin{equation}\label{eq:z0_uy}
				\norm{\hat z_\omega(0)}^2 = \int_{0}^{\infty} \inner{\hat u_\omega(-t)}{\hat y_\omega(t)}\, \d t
			\end{equation}
		Since the Fourier transform is unitary, it follows that 
		\begin{equation}\label{eq:z0}
				\norm{z(0)}^2_\calZ = \norm{\hat z(0)}^2_\calZ = \int_{\bbR^s} \norm{\hat z_\omega(0)}^2\, \d\omega
			\end{equation}
		and, similarly,
		\begin{multline}\label{eq:uy}
				\int_{0}^{\infty} \inner{u(-t)}{y(t)}_\calU\, \d t = \int_{\bbR^s} \int_{0}^{\infty} \inner{\hat u_\omega(-t)}{\hat y_\omega(t)}\, \d t \d\omega
			\end{multline}
		where we have used Fubini's theorem to interchange the order of integration. The right-hand sides of \eqref{eq:z0} and \eqref{eq:uy} are equal because of \eqref{eq:z0_uy}, hence the left-hand sides are also equal, which concludes the proof.
	\endproof
	
	\newcommand{\intr}{\int_{\bbR}}
	The right-hand side of \eqref{eq:ltsi_relaxation_passive} can be seen as the inner product between past input and future output. In fact, it is the quadratic functional obtained from the Hankel operator of the system, i.e., the operator that maps the past input to the future output. We formalize this in the next section. Before we do that, we illustrate the results of this section with the diffusion equation and then conclude with a couple of remarks.
	\begin{example}
		In Example~\ref{ex:diff_relaxation}, we saw that the diffusion equation \eqref{eq:diff} is an LTSI system that is internally of relaxation type. The state and output of this system are equal, hence, due to Theorem~\ref{thm:storage_io}, it is impedance passive with the well-known energy storage $E(t) = \int_\bbR y(t,x)^2\,\d x$ of the diffusion equation. In particular, \eqref{eq:ltsi_relaxation_passive} reduces to
		\begin{equation}
			\intr y(0,x)^2\, \d x = \int_0^\infty \intr u(-t,x)y(t,x)\, \d x\d t,
		\end{equation}
		which, due to Plancherel's theorem, is equivalent to
		\begin{equation}\label{eq:diff_y0omega_norm}
			\intr \abs{\hat y(0,\omega)}^2\, \d \omega = \int_0^\infty \intr \hat u(-t,\omega)\hat y(t,\omega)\a \, \d \omega\d t,
		\end{equation}
		Note that the future output $y(t,x)$, $t\geq0$, is the solution of the diffusion equation \eqref{eq:diff} for $u(t,x) = 0$, $t\geq0$, and initial condition $y(0,x)$. Therefore, in the frequency domain, we have that $\hat y(t,\omega) = \hat g(t,\omega)\hat y(0,\omega)$, where
		\begin{align}
			\hat y(0,\omega) = \int_{-\infty}^{0} \hat g(-t,\omega)\hat u(t,\omega)\, \d t.
		\end{align}
		Multiplying the latter by \new{$\hat y(0,\omega)\a$} and rearranging yields
		\begin{equation}\label{eq:diff_y0omega}
			\abs{\hat y(0,\omega)}^2 = \int_0^\infty  \hat u(-t,\omega)\hat y(t,\omega)\a  \, \d t,
		\end{equation}
		where we used the fact that $\hat g(t,\omega)$ is real and, thus,
		\begin{equation}
			\hat y(t,\omega)\a = \hat g(t,\omega) \hat y(0,\omega)\a.
		\end{equation}
		Integrating both sides of \eqref{eq:diff_y0omega} yields \eqref{eq:diff_y0omega_norm}, hence the energy storage of the diffusion equation can be determined \new{from} the past input, which confirms the result of Theorem~\ref{thm:storage_io}.
	\end{example}
	
	\begin{remark}\label{rem:relaxation_passivity}
		LTSI relaxation systems are not necessarily impedance passive. To see this, suppose that the LTSI system $\sys$ is of relaxation type, hence, due to Theorem~\ref{thm:relaxation_fourier}, the LTI system $\hat \sys_\omega$ is of relaxation type and, thus, passive for all $\omega\in\bbR^s$. This implies that there exist $\hat Q_\omega \in \bbC^{n\times n}$, $\hat Q_\omega \geq 0$, $\omega\in\bbR^s$, such that \eqref{eq:passivity_omega} holds, but it does not imply that $\sup_{\omega\in\bbR^s} \norm{\hat Q_\omega} < \infty$. It is argued in \cite{bamieh2002} that cases where such boundedness conditions are violated are artificial mathematical constructs rather than real physical examples. This is because at large frequencies, the dominant mechanism in physical systems is dissipation, i.e., the LTI system $\hat \sys_\omega$ becomes ``more stable'' as $\omega \to \infty$.
	\end{remark}
	\begin{remark}\label{rem:passivity_definitions}
		In some references \cite{arov2005, reis2008}, the definition of impedance passivity requires the storage to be strictly positive rather than merely nonnegative. In fact, it is sometimes assumed \cite{staffans2002} that the storage is the squared norm of the state. These definitions coincide with the definition in this paper if we further impose that $\calQ > 0$ such that $\norm{z}_\calQ = \sqrt{\inner{\calQ z}{z}}$ defines a norm on the state space $\calZ$. In any case, Theorem~\ref{thm:storage_io} shows that LTSI systems that are internally of relaxation type are impedance passive even under these modified definitions.
	\end{remark}

	\section{The Hankel operator and relaxation}\label{sec:hankel_relaxation}
	
	In this subsection, we show that relaxation can be characterized via the Hankel operator. In particular, we show that an exponentially stable LTSI system is of relaxation type if and only if its Hankel operator is self-adjoint and nonnegative. Consequently, the Hankel operator of an LTSI relaxation system is cyclically monotone and, thus, the gradient of a closed convex functional, which coincides with the (impedance) passive storage of Theorem~\ref{thm:storage_io}.
	
	The Hankel operator of the LTSI system $\sys$ is the linear operator $\calH:\D(\calH) \to L_2(\bbRp, \calY)$ given by
	\begin{equation}
		(\calH v)(t) = \intrp \calG(t+\tau)v(\tau)\, \d \tau,
	\end{equation}
	where $\D(\calH)\subset L_2(\bbRp,\calU)$. The Hankel operator maps past inputs to future outputs in the sense that the input given by $u(t) = v(-t)$ for $t\leq 0$ and $u(t) = 0$ for all $t >0$, results in the output given by $y(t) = (\calH v)(t)$ for all $t\geq 0$. The domain of $\calH$ is nontrivial only if the system is stable. Here, we assume that $\sys$ is \emph{exponentially stable}, i.e., there exists $M, \alpha >0$ such that $\norm{\exp(\calA t) }\leq Me^{-\alpha t}$ for all $t\in\bbRp$, see \cite[Section~4.1]{curtain2020} for details. Then, we obtain the following theorem.
	\begin{theorem}\label{thm:hankel_bounded}
		If the LTSI system $\sys$ is exponentially stable, then the Hankel operator $\calH$ extends to a bounded linear operator from $L_2(\bbRp,\calU)$ to $L_2(\bbRp,\calY)$.
	\end{theorem}
	\proof
		Note that $\exp(\calA(t+\tau)) = \exp(\calA t)\exp(\calA \tau)$ by definition of a semigroup, hence
		\begin{equation}
				(\calH v)(t) = \calC\exp(\calA t)\intrp\exp(\calA \tau)\calB v(\tau)\, \d \tau.
			\end{equation}
		Let the linear operator $\calB_\infty: L_2(\bbRp, \calU) \to\calZ$ be  given by
		\begin{equation}
				(\calB_\infty v)(t) = \intrp \exp(\calA \tau)\calB v(\tau)\, \d\tau
			\end{equation}
		and the linear operator $\calC_\infty: \calZ \to L_2(\bbRp,\calY)$  by
		\begin{equation}
				(\calC_\infty z)(t) = \calC\exp(\calA t)z
			\end{equation}
		The operators $\calB_\infty$ and $\calC_\infty$ are known as the extended controllability and observability maps, respectively, see \cite[Section~6.4]{curtain2020} for details. They are well-defined and bounded if $\sys$ is exponentially stable, hence $\calH = \calC_\infty \calB_\infty$ is a bounded linear operator from $L_2(\bbRp,\calU)$ to $L_2(\bbRp,\calY)$.
	\endproof
	
	\new{As shown in \cite{chaffey2023}, a stable} LTI system is of relaxation type if and only if its Hankel operator is self-adjoint and nonnegative (hence cyclically monotone). Here, we extend this result to exponentially stable LTSI systems.
	\begin{theorem}\label{thm:hankel_spd}
		Suppose that the LTSI system $\sys$ is exponentially stable. Then, $\sys$ is of relaxation type if and only if the Hankel operator $\calH\in\B(L_2(\bbRp,\calU))$ satisfies $\calH \geq 0$. 
	\end{theorem}
	\proof
		Note that $\calH\in\B(L_2(\bbRp,\calU))$ due to Theorem~\ref{thm:hankel_bounded} and the assumption that $\sys$ is exponentially stable. We first prove necessity. Suppose that $\sys$ is of relaxation type. It follows that $\hat\sys_\omega$ is of relaxation type for all $\omega\in\bbR^s$ due to Theorem~\ref{thm:relaxation_fourier}. Furthermore, due to \cite[Theorem~1]{bamieh2002}, $\hat \sys_\omega$ is stable for all $\omega\in\bbR^s$. Consequently, due to \cite[Theorem~5]{chaffey2023}, the Hankel operator $\hat\calH_\omega\in\B(L_2(\bbRp, \bbC^m))$ of $\hat \sys_\omega$ satisfies $\hat\calH_\omega\geq 0$ for all $\omega\in\bbR^s$. Note that
		\begin{equation}
				\inner{\hat \calH_\omega \hat v_\omega}{\hat w_\omega} = \intrp \intrp \inner{\hat g_\omega(t+\tau) \hat v_\omega(\tau)}{\hat w_\omega(t)}\, \d\tau\d t
			\end{equation}
		for all $\hat v_\omega,\hat w_\omega \in L_2(\bbRp, \bbC^m)$, and, similarly,
		\begin{equation}\label{eq:Hv_w}
				\inner{\calH v}{w} = \intrp \intrp \inner{\calG(t+\tau) v(\tau)}{w(t)}\, \d\tau\d t.
			\end{equation}
		for all $v,w\in L_2(\bbRp,\calU)$. Due to Plancherel's theorem,
		\begin{align*}
				\inner{\calG(t+\tau) v(\tau)}{w(t)} &=  \inner{\hat \calG(t+\tau) \hat v(\tau)}{\hat w(t)} \\
				&= \intrs  \inner{\hat g_\omega(t+\tau) \hat v_\omega(\tau)}{\hat w_\omega(t)}\, \d\omega
			\end{align*}
		for all $v,w\in L_2(\bbRp,\calU)$ with (spatial) Fourier transforms $\hat v,\hat w \in L_2(\bbRp, \calU)$.
		Substituting the latter in \eqref{eq:Hv_w} and changing the order of integration yields
		\begin{equation}\label{eq:Hv_w_fourier}
				\inner{\calH v}{w} = \intrs \inner{\hat\calH_\omega \hat v_\omega}{\hat w_\omega}\, \d\omega
			\end{equation}
		where $ \hat v_\omega,\hat w_\omega \in L_2(\bbRp, \bbC^m)$ for almost all $\omega\in\bbR^s$. Therefore, $\calH \geq 0$ because $\hat \calH_\omega \geq 0$ for all $\omega\in\bbR^s$.
		
		Next, we prove sufficiency. Suppose that $\calH \geq 0$. {Consider arbitrary $\bar v,\bar w\in L_2(\bbRp, \bbC^m)$. Let $\hat v_\omega = \bar v$ and $\hat w_\omega = \bar w$ for all $\omega\in\Omega$, and $\hat v_\omega = \hat w_\omega = 0$ for all $\omega\notin\Omega$, where $\Omega\subset\bbR^s$ is an arbitrary compact subset. Then, $v,w\in L_2(\bbRp,\calU)$ and
		\begin{equation}
			\int_\Omega \inner{\hat\calH_\omega \bar v}{\bar w} - \inner{ \bar v}{\hat\calH_\omega \bar w}\, \d\omega = 0
		\end{equation}
		due to \eqref{eq:Hv_w_fourier}  and the assumption that $\calH$ is self-adjoint. The latter holds for all $\bar v,\bar w\in L_2(\bbRp, \bbC^m)$ and compact $\Omega\subset\bbR^s$, hence }$\hat \calH_\omega$ is self-adjoint for almost all $\omega\in\bbR^s$. Similarly, $\hat \calH_\omega$ is nonnegative for almost all $\omega\in\bbR^s$. Due to \cite[Theorem~5]{chaffey2023}, this implies that $\hat \sys_\omega$ is of relaxation type and, thus, $\hat g_\omega$ is completely monotone for almost all $\omega\in\bbR^s$. Note that $\hat g_\omega(t)$ and its time derivatives are continuous in $\omega\in\bbR^s$ for each $t\in\bbRp$ because the symbols $\hat A$, $\hat B$ and $\hat C$ are assumed to be continuous. Therefore, by continuity, $\hat g_\omega$ is completely monotone for \emph{all} $\omega\in\bbR^s$, hence $\sys$ is of relaxation type due to Theorem~\ref{thm:relaxation_fourier}.
	\endproof
	
	As a consequence of Theorem~\ref{thm:hankel_spd}, the Hankel operator $\calH$ of an LTSI relaxation system $\sys$ is (maximally) cyclically monotone and, thus, the gradient of a closed convex functional \cite{rockafellar1966, rockafellar1970}. In fact, $\calH$ is the functional derivative of the convex quadratic functional $\frakH: L_2(\bbRp, \calU) \to \bbRp$ given by
	\begin{equation}
		\frakH(v) = \frac{1}{2} \inner{\calH v}{v},
	\end{equation}
	see \cite[Lemma~1]{chaffey2023}.  Theorem~\ref{thm:storage_io} shows that $\frakH$ is intimately related to the (impedance passive) storage of a relaxation system. In particular, we can write \eqref{eq:ltsi_relaxation_passive} as 
	\begin{equation}
		\calS(z(0)) = \inner{\calH \bar u}{\bar u} = 2\frakH(\bar u),
	\end{equation}
	where $\bar u(t) = u(-t)$ for all $t \in\bbRp$. The latter implies that the energy stored in a relaxation system is completely determined from the past input, i.e., its memory. In the LTI case \cite{chaffey2023}, the functional $\frakH$ defines a so-called \emph{intrinsic storage}, which is used to show passivity of relaxation systems from a purely input-output perspective. An analogous derivation is also possible in the LTSI case. As always, the main idea is that the LTSI system $\sys$ inherits the properties of the corresponding family of LTI systems $\hat\sys_\omega$, and vice versa.
	
	\section{Conclusion}\label{sec:conclusion}

	We extended the classical concept of relaxation from LTI to a class of LTSI systems. We showed that LTSI relaxation systems enjoy the same properties as LTI relaxation systems, i.e., they have completely monotone impulse responses, they are impedance passive, and their Hankel operators are cyclically monotone. As in the LTI case, LTSI relaxation systems reconcile the state-space concept of energy storage with the input-output concept of memory functional. Throughout this paper, we made use of the properties of the Fourier transform with respect to translation invariance to decouple the infinite-dimensional LTSI system into an infinite family of finite-dimensional LTI systems, so that properties of one can be characterized as properties of the other.
	
	We identify several topics for future research. It would be interesting to investigate whether the apparent discrepancy between internal and external relaxation can be resolved by imposing spatial properties on the impulse response, or, alternatively, by working in a Hardy space rather than $L_2$. It would also be interesting to consider LTSI relaxation from a purely geometric viewpoint \cite{vanderschaft2024}. In any case, we intend to extend the theory developed in this paper to nonlinear systems, in the same spirit as \cite{rodolphe2024}. In particular, we conjecture that spatio-temporal input-output operators that are \emph{myopic} in space \cite{sandberg1997} with completely monotone \emph{fading memory} \cite{boyd1985} can be universally approximated by an LTSI relaxation system composed with an absolutely monotone readout.
	
	\bibliographystyle{ieeetr}
	\bibliography{../ltsi_references}
\end{document}